\documentclass[11pt]{epiarticle}
\usepackage{epimath}
\setpapertype{A4}
\usepackage[english]{babel}
\usepackage[pdftex]{graphicx}

\title{A Well-Motivated Proof That Pi Is Irrational}
\titlemark{A Well-Motivated Proof That Pi Is Irrational}
\author{Timothy Y. Chow}
\date{March 29, 2024}
\journal{Journal -- (2024), ---} 

\begin{document}
\maketitle

\begin{prelims}
\def\abstractname{Abstract}
\abstract{Ivan Niven's succinct proof that $\pi$ is irrational is easy
to verify, but it begins with a magical formula that appears to come
out of nowhere, and whose origin remains mysterious even after one goes through
the proof. The goal of this expository paper is to describe a thought
process by which a mathematician might come up with the proof from scratch,
without having to be a genius. Compared to previous expositions of Niven's
proof, perhaps the main novelty in the present account is an explicit
appeal to the theory of orthogonal polynomials, which leads naturally to the
consideration of certain integrals whose relevance is otherwise not
immediately obvious.}

\end{prelims}

\section{Introduction}

That $\pi$ is irrational is something we have all known since childhood,
but curiously enough, most mathematicians have
either never seen a proof that $\pi$ is irrational,
or have worked through a proof
but have found it to be unmotivated and unmemorable.
This is not because of a lack of concise proofs;
Niven's famous proof (\cite{niven} or \cite[p.~276]{bbb}),
or a variant thereof \cite{bourbaki, hardy-wright,jeffreys},
occupies less than a page of text, and it is not difficult
to check that each step in the proof is correct.
A \emph{concise} proof,
however, is not the same as what Donald Newman~\cite{newman}
would call a \emph{natural} proof.
\begin{quote}
This term \dots\ is introduced to mean not having any ad hoc
constructions or \emph{brilliancies}. A ``natural'' proof, then, is one which
proves itself, one available to the ``common mathematician in the
streets.''
\end{quote}
Indeed, a verbose proof may be more natural than a concise proof,
if the concise proof fails to explain the origin of the underlying ideas.
For example, Niven's proof begins by writing down an integral that
appears to come out of nowhere.
Various authors \cite{jones,muller,zhou1,zhou2} have attempted to
motivate Niven's proof, but I have always been left with the feeling that
never in a million years could I have come up with the proof myself.

It was not until recently, after reading
Angell's lovely book~\cite{angell}, as well as the
answer by Kostya\_I~\cite{kostya} to a question that I
posted on MathOverflow, that a light bulb went off in
my head. The purpose of this paper is to give an account
of how Newman's ``mathematician in the streets'' might
discover a proof that $\pi$ is irrational
without having to be a genius.

In order to make this paper as accessible as possible,
we do not assume that the reader has any prior familiarity
with irrationality proofs.  Our discussion proceeds in several stages.
\begin{enumerate}
\item We explain the general philosophy behind irrationality proofs,
using Fourier's proof that $e$ is irrational as an example.
\item We give a proof that $e^r$ is irrational for positive integers~$r$
that minimizes ``brilliancies'';
ideally, readers will feel that they might have come up with this proof
themselves---at least if they were permitted to look up
``standard'' facts in ``standard'' references.
\item Using a clever trick, we simplify the proof, to enable the reader
to keep it in memory (or reconstruct it) without having to consult any
references.
\item Finally, we show that the proof that $e^r$ is irrational for
positive integers~$r$ can be
straightforwardly modified to yield Niven's proof that $\pi$ is irrational.
\end{enumerate}

\section{The fundamental theorem of transcendental number theory}

It is an old joke that the fundamental theorem of
transcendental number theory is that there is no integer strictly
between 0 and~1. Actually, this joke is half serious, because
many proofs of irrationality can be framed as follows.
\begin{enumerate}
\item Assume toward a contradiction that $\alpha$ is rational.
Write down a suitable equation involving~$\alpha$.
\item ``Scale up'' the equation by a multiple of the
denominator of~$\alpha$.
\item Deduce that $A=B$ where $A$ is an integer and $0<B<1$.
\item Apply the fundamental theorem of transcendental
number theory to deduce a contradiction.
\end{enumerate}
Fourier's proof~\cite{stainville} of the irrationality of~$e$ illustrates this schema
perfectly. For the purposes of this proof, we define the
function $e^x$ in terms of its Taylor series, so
\begin{equation}
\label{eq:e}
e = 1 + \frac{1}{1!} + \frac{1}{2!} + \frac{1}{3!} + \frac{1}{4!} + \cdots
\end{equation}
Assume toward a contradiction that $e = p/q$ for positive integers
$p$ and~$q$. Scale up both sides of Equation~\eqref{eq:e} by a
factor of~$q!$ and observe that some terms are integers;
let $B$ denote the sum of the remaining terms.
\begin{equation*}
\underbrace{\frac{q!p}{q}}_{\textstyle \in\mathbb{Z}} =
  \underbrace{q! + \frac{q!}{1!} + \frac{q!}{2!} + \cdots
   + \frac{q!}{q!}}_{\textstyle \in \mathbb{Z}} +
   \underbrace{\frac{q!}{(q+1)!} + \frac{q!}{(q+2)!} + \cdots}_{\textstyle B}
\end{equation*}
Since $B$ is the sum of positive terms, $0<B$.
It is straightforward to upper-bound $B$ with a geometric series,
which we can sum explicitly.
\begin{equation*}
B = \frac{1}{q+1} + \frac{1}{(q+1)(q+2)} + \cdots
  < \frac{1}{q+1} + \frac{1}{(q+1)^2} + \cdots
  = \frac{1}{q} \le 1.
\end{equation*}
This contradicts the fundamental theorem of
transcendental number theory; Q.E.D.

It is worth remarking explicitly that what makes this proof work is
the \emph{rapid convergence} of Equation~\eqref{eq:e}.
That is, no matter what $p/q$ we postulate $e$ to be, the ``remainder''
after scaling up by~$q!$ is so small that it lies between 0 and~1.

\section[The irrationality of \texorpdfstring{$e^r$}{n}]{The irrationality of \boldmath{$e^r$}}

Emboldened by our success, we might try to adapt Fourier's proof to
prove that $e^r$ is irrational when $r$ is a positive integer.
Let $e^r = p/q$ as before. Then
\begin{equation*}
\frac{p}{q} = 1 + \frac{r}{1!} + \frac{r^2}{2!} + \frac{r^3}{3!} + \cdots
\end{equation*}
What now? The most obvious attempt is to scale up by $q!$ as before.
\begin{equation*}
\underbrace{\frac{q! p}{q}}_{\textstyle \in\mathbb{Z}} =
\underbrace{q! + \frac{q! r}{1!}
   + \cdots + \frac{q! r^q}{q!}}_{\textstyle \in\mathbb{Z}} +
  \underbrace{\frac{q! r^{q+1}}{(q+1)!}
  + \frac{q! r^{q+2}}{(q+2)!}  + \cdots}_{\textstyle B}
\end{equation*}
But now we run into a difficulty.
If we try to upper-bound $B$ with a geometric series in the same
manner as before, then the bound we get is
\begin{equation*}
B < \sum_{n=1}^\infty \frac{r^{q+n}}{(q+1)^n }
  = \frac{r^{q+1}}{(q+1) - r}.
\end{equation*}
We get a contradiction if $r=1$,
but if $r>1$, then our upper bound on~$B$ is large
for large~$q$, so our argument fails.

It is quite possible that the above argument can be repaired, but
I do not know how to do so. For the purposes of the present discussion,
let us assume we get stuck at this point.
The Taylor series for~$e^x$ is not bringing us joy.
What alternatives exist?

\section{Orthogonal polynomials to the rescue}

We now come to a vital step in the proof.
Recall that we desire a rapidly converging approximation to~$e^x$.
Where might we find such a thing?
The idea is to expand
$e^x$ in a basis of \emph{orthogonal polynomials}.

Readers who are not very familiar with orthogonal polynomials will
have at least encountered \emph{Chebyshev polynomials of the
first kind} (if not by that name) in high-school trigonometry.
If we write $\cos n\theta$ as a polynomial in $\cos\theta$,
then the polynomials~$T_n$ that arise are precisely the
Chebyshev polynomials; e.g.,
\begin{align*}
T_2(x) = 2x^2 - 1 \quad\text{because} \quad
\cos 2\theta &= 2 (\cos\theta)^2 - 1,\\
T_3(x) = 4x^3 - 3x \quad\text{because} \quad
\cos 3\theta &= 4(\cos\theta)^3 - 3\cos \theta,\\
T_4(x) = 8x^4 - 8x^2 + 1 \quad\text{because} \quad
\cos 4\theta &= 8(\cos\theta)^4 - 8(\cos \theta)^2 + 1.
\end{align*}
The reason we say that the $T_n$ are orthogonal polynomials is that
they turn out to be orthogonal with respect to the following inner product:
\begin{equation*}
\langle f, g \rangle := \int_{-1}^1 f(x) g(x) \,
  \frac{dx}{\sqrt{1-x^2}}.
\end{equation*}
As in any inner product space, we can expand the vectors
(which in this case are functions of~$x$) in terms of
a basis of orthogonal vectors, and if the basis is a
``good'' one, then ``nice'' functions will be well-approximated.

Once we hit on the idea of computing the coefficients of~$e^x$
with respect to a basis of orthogonal polynomials, the question immediately
arises: Which set of orthogonal polynomials should we pick?
If we are not experts in the subject, we might turn to Wikipedia
or some standard reference such as the \emph{NIST Handbook}~\cite{nist}.
We quickly discover a bewildering variety
of orthogonal polynomials and associated inner products,
such as the following.
\begin{align*}
\text{Chebyshev polynomials:} \quad &
  \langle f, g \rangle := \int_{-1}^1 f(x) g(x) \, \frac{dx}{\sqrt{1-x^2}} \\
\text{Gegenbauer polynomials:} \quad &
   \langle f, g \rangle := \int_{-1}^1 f(x) g(x) (1-x^2)^\alpha \, dx \\
\text{Legendre polynomials:} \quad &
   \langle f, g \rangle := \int_{-1}^1 f(x) g(x) \, dx \\
\text{Laguerre polynomials:} \quad &
   \langle f, g \rangle := \int_0^\infty f(x) g(x) \,e^{-x}\, dx \\
\text{Hermite polynomials:} \quad &
   \langle f, g \rangle := \int_{-\infty}^\infty f(x) g(x)\, e^{-x^2}\, dx
\end{align*}
What to do?
In the absence of a strong a priori reason to pick one inner product
over another, it makes sense to begin with the simplest-looking
option---the inner product for Legendre polynomials---and switch horses
later if the first attempt does not work.
In this way, we arrive
at the idea of considering an expression of the
form $\int f(x)\,e^x\, dx$ where $f(x)$ is a polynomial.

Before we pull out the \emph{NIST Handbook} to look up
standard results on Legendre polynomials,
let us note that regardless of which polynomial $f(x)$ we choose,
evaluating $\int f(x)\,e^x\,dx$ is likely to require integrating by parts.
So let us see what we get if we follow our nose.
The parameter~$r$ needs to be introduced somehow,
so let us change the interval of integration from $[-1,1]$ to $[0,r]$.
\begin{align*}
\int_0^r f(x)\, e^x \, dx &= \biggl[f(x)\, e^x\biggr]_0^r - \int_0^r f'(x)\, e^x \, dx \\
&= \biggl[(f(x) - f'(x))\,e^x \biggr]_0^r + \int_0^r f''(x)\, e^x \, dx \\
&= \biggl[(f(x) - f'(x) + f''(x))\,e^x \biggr]_0^r - \int_0^r f'''(x)\, e^x \, dx,
 \quad \text{etc.}
\end{align*}
Therefore, if we define $F(x) := f(x) - f'(x) + f''(x) - f'''(x) + \cdots$
(there are no convergence problems, because $f(x)$ is a polynomial), then
\begin{equation}
\label{eq:intfxex}
\int_0^r f(x) \, e^x \, dx = F(r)\, e^r - F(0).
\end{equation}
If the fundamental theorem of transcendental number theory is burned
into our minds, Equation~\eqref{eq:intfxex} should make us perk up.
Suppose toward a contradiction that $e^r = p/q$.
Can we choose $f(x)$ so that, after scaling up by~$q$,
the right-hand side is an integer,
but the left-hand side is strictly between 0 and~1?
There seems to be a tension; the simplest way to ensure that
the right-hand side of Equation~\eqref{eq:intfxex} is an integer
(after multiplying by~$q$) is to 
require $f(x)$ to have integer coefficients,
because then $F(r)$ and $F(0)$ will be integers.
But if $f(x)$ has integer coefficients,
then there would appear to be no reason for the integral
on the left-hand side to be small.

\section{Legendre polynomials}
\label{sec:legendre}

It is time to look up the theory of Legendre polynomials,
to see if it can help us.
The usual definition of Legendre polynomials assumes
an interval of integration of $[-1,1]$;
one definition~\cite[18.5.8]{nist} is
\begin{equation*}
P_n(x) = \frac{1}{2^n} \sum_{\ell=0}^n
  \left(\genfrac{}{}{0pt}{}{n}{\ell}\right)^2
  (x-1)^{n-\ell}(x+1)^\ell.
\end{equation*}
Recall that we have changed the interval of integration to $[0,r]$,
so we need to make an affine change of variables. Define
\begin{equation}
\label{eq:tildelegendre}
\tilde P_n(x) := P_n(2x/r - 1) = \frac{1}{r^n} \sum_{\ell=0}^n
  \left(\genfrac{}{}{0pt}{}{n}{\ell}\right)^2
  x^\ell (x-r)^{n-\ell}.
\end{equation}
Now, the $n$th coefficient in the Legendre polynomial expansion
of a function such as~$e^x$ is (up to a normalizing factor that
we will ignore for now) obtained by taking the
inner product of~$e^x$ with the $n$th Legendre polynomial:
\begin{equation*}
\int_0^r \tilde P_n(x) \, e^x \, dx.
\end{equation*}
As we alluded to earlier, it is natural
to guess that Legendre polynomials form a
``good'' basis in the sense that these coefficients shrink
rapidly as $n$ increases.
Aha! We want the integral in Equation~\eqref{eq:intfxex}
to be small, so maybe setting $f(x) := \tilde P_n(x)$ is the key?
But remember that we also want $f(x)$ to have integer coefficients.
Equation~\eqref{eq:tildelegendre} implies that
$r^n \tilde P_n(x)$ has integer coefficients,
so let us try setting $f(x) := r^n \tilde P_n(x)$.
Then Equation~\eqref{eq:intfxex} becomes (after multiplying by~$q$)
\begin{equation*}
r^n q \int_0^r \tilde P_n(x) \, e^x \, dx
 = F(r)\cdot p - F(0)\cdot q,
\end{equation*}
where $F(x) = f(x) - f'(x) + f''(x) - \cdots\,$.
The right-hand side is an integer,
so to obtain our desired contradiction, it would be enough if
\begin{equation}
\label{eq:innerproduct}
r^n q \int_0^r \tilde P_n(x) \, e^x \, dx
\end{equation}
were nonzero, with absolute value less than~$1$ for some~$n$.
But is this the case?

It turns out that the answer is yes!
The proof is a relatively routine but
somewhat tedious calculation that we will present in a moment.
We therefore have:

\bigskip
\noindent\textbf{Theorem 1.} 
\emph{If $r$ is a positive integer, then $e^r$ is irrational.}
\bigskip

Let us pause a moment to take stock of the situation.
We claim that we have managed to find a proof that $e^r$ is irrational
(for positive integers~$r$) with a minimum of flashes of brilliance.
Once we think to use orthogonal
polynomials---and in particular Legendre polynomials---to
give us a good approximation to~$e^x$, we are led to
consider an integral (over a finite interval) of the form
$\int f(x)\,e^x\, dx$ where $f(x)$ is a
(suitably normalized) Legendre polynomial.
On the one hand, we can calculate directly that the integral,
even after multiplication by the denominator of~$e^r$,
shrinks rapidly. On the other hand, we can use integration
by parts to show that its value is an integer,
contradicting the fundamental theorem of
transcendental number theory.

\begin{proof}[First proof of Theorem~1]
The following argument is not needed for anything else in this paper,
so the reader who is unfamiliar with orthogonal polynomials is
encouraged to skip ahead to Section~\ref{sec:shorter} and
return here later.

As we said above, we just need to prove that
\eqref{eq:innerproduct} is nonzero and has
absolute value less than~1 for some~$n$.
The reason \eqref{eq:innerproduct} is nonzero for infinitely many~$n$
is that Legendre polynomials form an orthogonal basis
for the Hilbert space of square-integrable functions on a compact interval,
so having zero coefficients for all sufficiently large~$n$
would imply that $e^x$ is a \textit{finite} sum of polynomials,
which of course is absurd.

The polynomial $\tilde P_n(x)$ is orthogonal to all
polynomials of degree less than~$n$,
so the value of the integral~\eqref{eq:innerproduct}
remains unchanged if we replace $e^x$ by $e_n(x)$,
which we define to be $e^x$ minus the terms in its
Taylor series of degree less than~$n$. The
Cauchy--Bunyakovsky--Schwarz inequality implies that
\begin{align*}
\left|\int_0^r \tilde P_n(x) \, e^x \, dx\right|
  &= \left|\int_0^r \tilde P_n(x) \, e_n(x) \, dx\right| \\
  &\le \int_0^r \left| \tilde P_n(x) \, e_n(x)\right| \, dx \\
  &\le ||\tilde P_n(x)||_2 \cdot ||e_n(x)||_2,
\end{align*}
where $||\cdot||_2$ denotes the $L^2[0,r]$ norm.
Using the substitution $y=2x/r - 1$,
\begin{equation*}
||\tilde P_n(x)||^2_2 = \int_0^r \tilde P_n(x)^2\, dx
 = \frac{r}{2} \int_{-1}^1 P_n(y)^2\, dy = \frac{r}{2n+1},
\end{equation*}
where the last equality above is a standard
fact~\cite[Table 18.3.1]{nist} about the
squared norm of the Legendre polynomial.
As for $e_n(x)$, for $0\le x \le r$ and large~$n$,
\begin{align*}
e_n(x) &= \frac{x^{n}}{n!} + \frac{x^{n+1}}{(n+1)!} + \cdots
     \le \frac{x^{n}}{n!} \biggl(1 + \frac{x}{n+1}
     + \frac{x^2}{(n+1)^2} + \cdots \biggr) \\
   &= \frac{x^{n}}{n!(1 - x/(n+1))} \le \frac{2x^{n}}{n!}.
\end{align*}
Therefore
\begin{align*}
||e_n(x)||_2 &= \sqrt{\int_0^r e_n(x)^2\, dx} 
             \le \frac{2}{n!} \sqrt{\int_0^r x^{2n}\, dx}\\
    &= \frac{2r^{(2n+1)/2}}{n!\sqrt{2n+1}}.
\end{align*}
Putting this all together yields the following upper
bound on the absolute value of~\eqref{eq:innerproduct}:
\begin{equation*}
r^n q \cdot \sqrt{\frac{r}{2n+1}} \cdot 
    \frac{2r^{(2n+1)/2}}{n!\sqrt{2n+1}},
\end{equation*}
which for fixed~$r$ and large~$n$ is less than~$1$,
because exponential functions grow much more slowly
than factorials.\qed
\end{proof}

\section{A shorter and more elementary proof}
\label{sec:shorter}

Readers who are satisfied that the proof we have just presented is
well-motivated may nevertheless be dissatisfied
that it requires us to know (or re-derive) various
facts about Legendre polynomials, and that the
calculations are a bit tedious.
The conceptual outline of the proof is not hard to remember,
and it motivates why we want to consider an integral
of the form $\int f(x) \, e^x\, dx$,
but the details are not so memorable.

In this section, we rectify these defects
with an alternative version of the proof
that is simpler and requires no ``theory.''
The price we pay is that we will need to invoke a clever idea.

Let us revisit Equation~\eqref{eq:intfxex}.
Assuming $e^r = p/q$ and $f(x)$ is a polynomial,
\begin{equation}
\label{eq:intfxex2}
q \int_0^r f(x) \, e^x \, dx = F(r)\cdot p - F(0) \cdot q.
\end{equation}
where $F(x) := f(x) - f'(x) + f''(x) - f'''(x) + \cdots\,$.
Recall that we want to choose $f(x)$ so that
the left-hand side of Equation~\eqref{eq:intfxex2}
is small (between 0 and~1),
and at the same time the right-hand side is an integer.
Can we figure out from first principles how to choose a suitable $f(x)$?

Here is the clever idea:
the definition of~$F(x)$ involves high-order derivatives~$f^{(n)}(x)$,
and \textit{repeated differentiation of
a polynomial produces a factorial-like coefficient out in front!}
That suggests that if $f(x)$ has integer coefficients,
then not only will $f^{(n)}(x)$ have integer coefficients,
but it will still have integer coefficients even after dividing by~$n!$.
Indeed:

\bigskip
\noindent\textbf{Lemma 1.}
\emph{If $f(x)$ is a polynomial with integer coefficients, then
for any nonnegative integer~$n$,
$f^{(n)}(x)/n!$ has integer coefficients.}

\begin{proof}
It suffices to prove Lemma~1 in the case $f(x) = x^m$.
The $n$th derivative of~$x^m$
is $m(m-1)\cdots(m-n+1)x^{m-n}$, and dividing 
the coefficient by~$n!$ yields
the binomial coefficient~$\binom{m}{n}$, which is an integer
for any nonnegative integer~$m$.\qed
\end{proof}

Lemma~1 is promising, because dividing by~$n!$ should make
the left-hand side of Equation~\eqref{eq:intfxex2} small.
But there is still a difficulty; if $k\ge n$, then we can divide
$f^{(k)}(x)$ by $n!$ and still get integer coefficients,
but what do we do about those pesky lower-order derivatives ($k < n$)
that show up in~$F(x)$?
We solve this problem by taking
advantage of our freedom to choose~$f(x)$.
If $f(x)$ vanishes to order~$n$ at some point~$x=c$,
then $f^{(k)}(c)$ vanishes for $k<n$.
Since the right-hand side of Equation~\eqref{eq:intfxex2} involves
$F(0)$ and $F(r)$, we are motivated to consider the
polynomial $x^n(r-x)^n$ (or the polynomial $x^n(x-r)^n$,
but we will see soon why $x^n(r-x)^n$ better serves our
purposes), which vanishes to order~$n$ at
both $x=0$ and $x=r$.

\bigskip
\noindent\textbf{Corollary 1.}
\emph{For integers $r$ and~$n$ with $n\ge 0$,
let $f(x) = x^n(r-x)^n$, and $F(x) = f(x) - f'(x) + f''(x) - \cdots\,$.
Then $F(0)/n!$ and $F(r)/n!$ are integers.}

\begin{proof}
Since $f(x)$ vanishes to order~$n$ at $x=0$ and $x=r$,
$f^{(k)}(0) = f^{(k)}(r) = 0$ if $k<n$. If $k\ge n$, then
Lemma~1 tells us that $f^{(k)}(x)/n!$ has integer coefficients.
Thus each summand of $F(0)/n!$ and $F(r)/n!$ is either zero or
an evaluation of a polynomial with integer coefficients
at an integer value.\qed
\end{proof}

%
%
%
%

We are now in a position to give a second proof of Theorem~1
that does not require any knowledge of orthogonal polynomials.

\begin{proof}[Second proof of Theorem~1]
Let $f(x) := x^n(r-x)^n$ for some $n$
(whose value will be chosen later),
and define $F(x) := f(x) - f'(x) + f''(x) - f'''(x) + \cdots\,$.
Assume toward a contradiction that
$e^r = p/q$ for positive integers $p$ and~$q$.
Equation~\eqref{eq:intfxex2} implies
\begin{equation}
\label{eq:theorem1}
\frac{q}{n!}\int_0^{r} x^n(r-x)^n \, e^x \, dx
  = \frac{F(r)p}{n!} - \frac{F(0)q}{n!}.
\end{equation}
Corollary~1 tells us that the right-hand side of
Equation~\eqref{eq:theorem1} is an integer.
On the left-hand side,
the interval of integration has nonzero length because $r>0$,
and the integrand is strictly positive between $0$ and $r$
(this is where our choice of $x^n(r-x)^n$ instead of
$x^n(x-r)^n$ makes a difference),
so the left-hand side is strictly positive.
But we can also upper-bound the integral
by multiplying the length of the interval of integration
by a term-by-term upper bound on the factors of the integrand:
\begin{equation*}
\int_0^{r} x^n(r-x)^n \, e^x \, dx \le r \cdot r^n \cdot r^n \cdot e^r.
\end{equation*}
Because $n!$ grows super-exponentially fast, the left-hand
side of Equation~\eqref{eq:theorem1} is less than~$1$ for
all sufficiently large~$n$, contradicting the fundamental
theorem of transcendental number theory.\qed
\end{proof}

We remark in passing that the two proofs of Theorem~1 are
not as different as they might seem at first glance,
because the Rodrigues formula~\cite[18.5.5]{nist}
for Legendre polynomials implies
\begin{equation*}
\tilde P_n(x) = \frac{1}{n!r^n} \frac{d^n}{dx^n} (x^n(x-r)^n),
\end{equation*}
so $\tilde P_n(x)$ is in fact closely related to the 
polynomial $f(x) = x^n(r-x)^n$.

\section[On to the irrationality of \texorpdfstring{$\pi$}{n}]{On to the irrationality of \boldmath{$\pi$}}

We have worked rather hard, and have not even mentioned $\pi$ yet.
But as we promised earlier, a small modification of the ideas we
have already seen will prove the irrationality of~$\pi$.

First, though, we should at least briefly address the basic question,
just what \emph{is} $\pi$ anyway?
The least controversial definition is that $\pi$ is the limit
of the perimeter of a regular $n$-gon of unit diameter
as $n$ goes to infinity. Unfortunately, this geometric definition is
awkward to work with. Instead, what we do is to define
\begin{equation*}
\sin x := x - \frac{x^3}{3!} + \frac{x^5}{5!} - \frac{x^7}{7!} + \cdots
\end{equation*}
and we define $\pi$ to be the smallest positive number such that
$\sin \pi = 0$. The advantage of this approach is that it allows
us to apply the tools of calculus to study~$\pi$.
Now, one could object that it is far from obvious that this
definition of~$\pi$ coincides with the goemetric one.
This objection is valid, but we take the view that the equivalence
of the two definitions is a standard fact that the
``common mathematician in the streets'' would be familiar with.

Theorem~1 says that $e^r$ is irrational for all positive integers~$r$,
but in fact it immediately implies
that $e^r$ is irrational for all nonzero
rational~$r$, because if $e^r$ is irrational then
so is $e^{-r} = 1/e^r$, and if $e^{r/a}$ is rational then
so is $(e^{r/a})^a = e^r$.
We can restate this result in the following form: if $r$ is nonzero
and $e^r$ is rational, then $r$ is irrational.
To prove that $\pi$ is irrational, 
\emph{we run through almost the same argument
with the function $\sin x$ instead of
the function $e^x$, and we take $r=\pi$}.
Now $\sin r = \sin\pi = 0$ is manifestly rational,
and we will argue that this implies that $\pi$ is irrational.

There are only two adjustments that need to be made in the proof.
First, replacing $e^x$ with $\sin x$ changes the integration by parts formula slightly.
\begin{align*}
\int_0^\pi f(x)\,\sin x\, dx
  &= \biggl[-f(x) \,\cos x\biggr]_0^\pi + \int_0^\pi f'(x)\,\cos x\, dx \\
  &= f(\pi) + f(0) + \biggl[f'(x)\, \sin x\biggr]_0^\pi
        - \int_0^\pi f''(x)\,\sin x\, dx\\
  &= f(\pi) + f(0) - \int_0^\pi f''(x)\,\sin x\, dx.
\end{align*}
So if we define $F(x) := f(x) - f''(x) + f''''(x) - f''''''(x) + \cdots$ then
\begin{equation}
\label{eq:intfxsinx}
\int_0^\pi f(x) \, \sin x\, dx = F(\pi) + F(0).
\end{equation}

The second adjustment arises because
we will be assuming toward a contradiction that
$\pi = a/b$, and since $\pi$ is playing the role that $r$ was
playing, we need a version of Corollary~1
in which $r-x$ is replaced with $a-bx$.

\bigskip
\noindent\textbf{Corollary 2.}
\emph{For integers $a, b, n$ with $b\ne 0$ and $n\ge 0$,
let $f(x) = x^n(a-bx)^n$, and $F(x) = f(x) - f''(x) + f''''(x) - \cdots\,$.
Then $F(0)/n!$ and $b^nF(a/b)/n!$ are integers.}

\bigskip
\noindent\emph{Proof.}
As in the proof of Corollary~1,
since $f(x)$ vanishes to order~$n$ at $x=0$ and $x=a/b$,
$f^{(k)}(0) = f^{(k)}(a/b) = 0$ if $k<n$. If $k\ge n$, then
Lemma~1 tells us that $f^{(k)}(x)/n!$ has integer coefficients.
Each summand of $F(0)/n!$ is either zero or the constant term
of a polynomial with integer coefficients, so $F(0)/n!$ is an integer.
As for $F(a/b)/n!$, it is the result of evaluating a degree-$n$ polynomial in~$x$
with integer coefficients at $x=a/b$, so multiplying it by~$b^n$
yields an integer. Q.E.D.

\bigskip
\noindent\textbf{Remark.} In fact, something even stronger is true;
the hypotheses of Corollary~2 imply that $F(a/b)/n!$ is already
an integer, because $n$-fold differentiation of $f(x)$ produces
suitable powers of~$b$ via the chain rule. But we do not need this
stronger result.

\bigskip
\noindent\textbf{Theorem 2.} 
\emph{$\pi$ is irrational.}

\begin{proof}
Assume toward a contradiction that
$\pi = a/b$ for positive integers $a$ and~$b$.
We define the polynomial $f(x) := x^n(a-bx)^n$,
and set $F(x) := f(x) - f''(x) + f''''(x) - f''''''(x) + \cdots\,$.
By Equation~\eqref{eq:intfxsinx},
\begin{equation*}
\int_0^\pi f(x) \, \sin x\, dx = F(\pi) + F(0).
\end{equation*}
If we multiply both sides of the equation by~$b^n\!/n!$, then we obtain
\begin{equation}
\label{eq:theorem2}
\frac{b^n}{n!}\int_0^{a/b} x^n(a-bx)^n \, \sin x \, dx
  = \frac{b^n F(a/b)}{n!} + \frac{b^nF(0)}{n!}.
\end{equation}
Corollary~2 implies that the right-hand side of
Equation~\eqref{eq:theorem2} is an integer.
On the left-hand side, the integrand is strictly positive
between $0$ and $a/b$ (this is where our decision to
use the function $\sin x$, rather than $\cos x$ or some
other phase shift, makes a difference), and the interval of integration is
nonzero because $\pi>0$, so the left-hand side is strictly positive.
But we can also upper-bound the integral
by multiplying the length of the interval of integration
by a term-by-term upper bound on the factors of the integrand:
\begin{equation*}
\int_0^{a/b} x^n(a-bx)^n \, \sin x \, dx \le
  \frac{a}{b}\biggl(\frac{a}{b}\biggr)^na^n.
\end{equation*}
Because $n!$ grows super-exponentially fast,
the left-hand side of Equation~\eqref{eq:theorem2} is less than~$1$ for
all sufficiently large~$n$, contradicting the fundamental
theorem of transcendental number theory.\qed
\end{proof}

We should remark that it is not hard to make similar adjustments
to the first proof of Theorem~1,
in Section~\ref{sec:legendre}, to prove the irrationality of~$\pi$
(see~\cite{kostya} for details).
Some readers may prefer this proof, if the ``clever idea''
of Section~\ref{sec:shorter} seems too ad hoc.

\section{Cheat sheet}

Despite the brevity of Niven's proof,
I was never previously able to retain it in my
long-term memory, because I did not understand
where all the ingredients came from.
But now, I am confident that I could reconstruct it
without having to look anything up. I have found that the following
``cheat sheet'' is all that I need to remember,
and I offer it to the reader as a mnemonic aid.

\begin{enumerate}
\item Motivated by the theory of orthogonal polynomials,
we consider an expression of the form $\int f(x)\,\sin x\,dx$
where $f(x)$ is a suitably chosen polynomial.
\item Integration by parts implies that if
$F(x) := f(x) - f''(x) + f''''(x) - \cdots$ then
\begin{equation*}
\int_0^\pi f(x) \, \sin x \, dx = F(\pi) + F(0).
\end{equation*}
\item If $\pi = a/b$, then we want the right-hand side
(after suitable scaling) to be an integer,
and the left-hand side to be between 0 and~1.
We do this by choosing $f$ to vanish to high order~$n$ at 0 and~$a/b$,
and then multiplying both sides by~$b^n\!/n!$.
The denominator of~$n!$ makes the LHS small.
The high-order vanishing implies that the RHS is an integer,
because low-order derivatives of~$f$ vanish at $0$ and~$\pi$,
and high-order derivatives of~$f$ have coefficients that
are divisible by~$n!$.
\end{enumerate}

\section{Historical remarks and acknowledgments}

As we mentioned earlier, the irrationality proof of~$\pi$
that we present here is due to Niven~\cite{niven},
although the ideas can be traced back to Hermite~\cite{hermite}
and \cite[pp.~162--193]{bbb}.
Our exposition, except for the connection to orthogonal
polynomials, borrows heavily from Angell~\cite{angell}.
The proof using Legendre polynomials is adapted from the
answer by Kostya\_I~\cite{kostya} to a question that I
posted on MathOverflow, asking about ways to motivate
the proof that $\pi$ is irrational.
Kostya\_I's proof was a revelation to me,
and a major inspiration for my writing this article in the first place.

The educated reader will notice that this paper
conspicuously fails to mention continued fractions,
which play an important role in the study
of rational approximations of irrational numbers.
Continued fractions were used by Lambert~\cite[pp.~141--146]{bbb}
in the first ever proof that $\pi$ is irrational,
and also show up in Hermite's work in the form of what
we now call \emph{Pad\'e approximants} (although Hermite
came before Pad\'e). Pad\'e approximants in turn are closely
related to orthogonal polynomials, so the relevance of
orthogonal polynomials to the irrationality of~$\pi$
is not unexpected, but we can perhaps claim some novelty in
the way we have emphasized their usefulness for expository purposes.
For more discussion of the relationship between
Equation~\eqref{eq:intfxex} and Pad\'e approximants
and the continued fraction
expansion of~$e$, see Cohn's \emph{Monthly} note~\cite{cohn}.

We have emphasized the parallel between the irrationality of~$e^r$
and the irrationality of~$\pi$, but we should point out that
the seemingly innocuous contraposition when switching from
``$r$~rational implies $e^r$ irrational'' to
``$\sin\pi$ rational implies $\pi$ irrational''
means that the proof technique yields explicit rational
approximations to~$e^r$, but does \emph{not} yield explicit
rational approximations to~$\pi$. Finding high-quality
explicit rational approximations to~$\pi$ is a topic
that is beyond the scope of this article.

Finally, we thank Nicholas Nguyen, John Zuehlke, and
especially Peter Farbman for carefully reading drafts of this paper;
their input has significantly improved the exposition.

\vfill\eject

\end{document}